\documentstyle[11pt]{article}
\begin{document}
\setlength{\baselineskip}{17pt}
\begin{center}
{\Large Local stability and Lyapunov functionals \\
for $n$-dimensional quasipolynomial conservative systems }
\end{center}

\mbox{}

\begin{center}
{\sc Benito Hern\'{a}ndez--Bermejo$^{\: 1}$ \hspace{2cm} 
V\'{\i}ctor Fair\'{e}n}
\end{center}

\mbox{}

\mbox{}

\noindent {\em Departamento de F\'{\i}sica Matem\'{a}tica y Fluidos, 
Universidad Nacional de Educaci\'{o}n a Distancia. Senda del Rey S/N, 28040 
Madrid, Spain.}

\mbox{}

\mbox{}

\mbox{}

\begin{center} 
{\sc Abstract}
\end{center}

We present a method for determining the local stability of equilibrium 
points of conservative generalizations of the Lotka-Volterra equations. These 
generalizations incorporate both an arbitrary number of species ---including 
odd-dimensional systems--- and nonlinearities of arbitrarily high order in 
the interspecific interaction terms. The method combines a 
reformulation of the equations in terms of a Poisson structure and the 
construction of their Lyapunov functionals via the energy-Casimir method. 
These Lyapunov functionals are a generalization of those traditionally known 
for Lotka-Volterra systems. Examples are given. 

\mbox{}

\noindent {\bf Running title:} Stability of nonlinear conservative systems.

\mbox{}


\mbox{}

\mbox{}

\mbox{}

\mbox{}

\mbox{}

\noindent $^1$ Corresponding author. E-mail: bhernand@apphys.uned.es \newline
Telephone: (+ 34 91) 398 72 19. \hspace{3mm} Fax: (+ 34 91) 398 66 97.

\pagebreak
\begin{center} 
{\sc 1. Introduction}
\end{center}

Consider the following Lotka-Volterra system 
[\ref{lot1},\ref{v1}]
\begin{equation}
   \label{lv0}
   \dot{x}_i = x_i \left( \lambda _i + \sum _{j=1}^n A_{ij} x_j \right) 
   , \;\:\;\: i = 1, \ldots ,n
\end{equation}
which is assumed to have a unique equilibrium point, ${\bf x}_0 \in 
\mbox{int}\{ I\!\!R_+^n\}$. One of the most relevant results about its  
stability is well summarized in a theorem originally enunciated by 
Kerner [\ref{ker1}], and later generalized by many different authors 
[\ref{goh1},\ref{goh2},\ref{har1},\ref{hsu1},\ref{kri1},\ref{ltj1},\ref{lyt1},\ref{red1},\ref{ta1}-\ref{tat2}] 
(see also [\ref{hs1},\ref{log1},\ref{tak1}] for detailed reviews of the 
subject). The result makes use of the well-known Lyapunov functional 
\begin{equation}
    \label{hlv0}
    V({\bf x}) = \sum _{i=1}^{n} d_i \left( x_i-x_{0i} -
      x_{0i} \ln \frac{x_i}{x_{0i}} \right)
\end{equation}
The time derivative of (\ref{hlv0}) along the trajectories of (\ref{lv0}) is
\begin{equation}
    \label{hp0}
    \dot{V}({\bf x}) = \frac{1}{2} ({\bf x} - {\bf x}_0)^T 
      (D \cdot A + A^T \cdot D) ({\bf x} - {\bf x}_0)
\end{equation}
where $D=$ diag$(d_1, \ldots ,d_n)$. Thus, it can be 
stated that if there exists a positive definite diagonal matrix $D$ such that 
$D \cdot A + A^T \cdot D$ is negative definite, ${\bf x}_0$ is Lyapunov 
asymptotically stable. Moreover if, instead, $D \cdot A + A^T \cdot D$ is 
negative semi-definite, then ${\bf x}_0$ is Lyapunov semi-stable and 
every solution of (\ref{lv0}) in int$\{ I\!\!R_+^n \}$ tends to the maximal 
invariant set $M$ contained in the set (see [\ref{las0},\ref{las1},\ref{lyt2}] 
and references therein)
\begin{equation}
   E = \left\{ {\bf x} \in \mbox{int} \{ I\!\!R_+^n \} \; / \; 
      ({\bf x} - {\bf x}_0)^T (D \cdot A + A^T \cdot D) ({\bf x} - 
      {\bf x}_0) = 0 \right\}
\end{equation}

Every one of the two previous alternatives encompasses the already classical 
community models, respectively: The so-called Lotka-Volterra dissipative and 
conservative systems [\ref{log1}]. In particular, Lotka-Volterra 
conservativeness implies that (\ref{hlv0}) is a constant of motion, thus 
making conservative systems formally amenable to analysis by standard 
theoretical mechanics methods [\ref{gas1},\ref{ker3},\ref{ker4}] and 
statistical mechanics considerations [\ref{gmm1},\ref{ker1},\ref{ker2}]. On 
this respect, the Hamiltonization of classical Lotka-Volterra conservative 
systems [\ref{ker3},\ref{ker4}] proceeds by defining the canonical variables, 
$z_i$, as linear transforms of new dependent variables $y_i= \ln (x_i/x_{0i})$. 
If the Hamiltonian, $H$, is simply identified with the functional in 
(\ref{hlv0}) appropriately rewritten in the canonical variables, $z_i$, the 
conservative Lotka-Volterra equations adopt the familiar 
symplectic form 
\begin{equation}
   \label{slp0}
   \dot{\bf z} = S \cdot \nabla H({\bf z})
\end{equation}
where $S$ is the classical symplectic matrix [\ref{gold}]. 

Unfortunately, this constructive procedure cannot be carried out beyond  
the class of even-dimensional classical conservative Lotka-Volterra systems. 
This made the Hamiltonian description of rather limited use until Nutku 
[\ref{nut1}] and Plank [\ref{pla1}] suggested reconsidering it under the 
light of the more general Poisson structure representation (see 
[\ref{olv1}] for an overview; see also references therein). Poisson systems 
constitute a natural extension of classical Hamiltonian dynamical systems, 
but have the advantage of embracing odd-dimensional flows as well. In the 
Poisson context, no prior transformation on the variables is necessary, and 
the conservative Lotka-Volterra equations can be put into 
Poisson form in terms of the original variables
\begin{equation}
   \label{p0}
   \dot{\bf x} = {\cal J} \cdot \nabla H({\bf x})
\end{equation}
where the elements of the structure matrix ${\cal J}$ are defined as 
$J_{ij} = K_{ij}x_ix_j$, $K$ being a skew-symmetric matrix, and $H$ is the 
classical Volterra's constant of motion (\ref{hlv0}). 

In fact, form (\ref{p0}) happens to be suitable for embracing a higher 
number of families of conservative systems than those of type (\ref{lv0}), as 
stated in the following result (see [\ref{byv5}]):

\mbox{}

{\sc Theorem 1.} [\ref{byv5}] \hspace{0.2cm} {\em Let us consider a 
differential system defined in the positive orthant, of the form 
\begin{equation}
   \dot{x}_i = x_i \left(\lambda _{i} + \sum_{j=1}^{m}A_{ij}\prod_{k=
      1}^{n}x_k^{B_{jk}} \right) , \;\:\;\: i = 1, \ldots ,n \; , \;\: m \geq n
   \label{glv}
\end{equation}
such that rank($B$) $=n$ and 
\begin{equation}
\label{kdl}
   \lambda = K \cdot L \; , \;\:\;\; A = K \cdot B^T \cdot D \; ,
\end{equation}
with $K$, $L$ and $D$ matrices of real entries, where $K$ is $n \times n$ and 
skew-symmetric; $L$ is $n \times 1$; and $D$ is $m \times m$, diagonal and of 
maximal rank. Then the system has a constant of motion of the form:
\begin{equation}
\label{H}
   H= \sum_{i=1}^m D_{ii} \prod_{k=1}^n x_k^{B_{ik}} + \sum_{j=1}^n L_j 
   \ln x_j 
\end{equation}
Moreover, the system is Poisson with Hamiltonian $H$ and a structure matrix 
${\cal J}$ with entries $J_{ij}= K_{ij}x_ix_j$.}

\mbox{}

Note that systems (\ref{glv}) appear when we combine a quadratic structure 
matrix (first identified by Plank [\ref{pla1}]) together with Hamiltonian 
(\ref{H}), which is a generalization of Volterra's constant of motion 
(\ref{hlv0}). Important dynamical features of certain particular cases of 
such systems have recently deserved detailed attention in the literature 
[\ref{pla2}]. In what follows, we shall denote systems described by Theorem 1 
as quasi\-polynomial of Poisson form, or QPP in brief. QPP systems (\ref{glv}) 
include the conservative Lotka-Volterra equations as a particular 
case when $m=n$, $B$ is the identity matrix, the dimension is even and $A$ is 
invertible. In such a case, Hamiltonian (\ref{H}) also reduces to Volterra's 
first integral, as it can be easily verified. 

The purpose of the present article is to investigate under which conditions 
the equilibrium points of the QPP systems are stable and compare the 
resulting generalization with what is known for conservative Lotka-Volterra 
models (\ref{lv0}). In particular, we shall also carry out a generalization 
of the corresponding Lyapunov functionals (\ref{hlv0}). In this way, we 
shall complete a treatment that simultaneously embraces arbitrary-dimensional 
systems and also arbitrary nonlinearities in the flow. 

The construction 
of suitable Lyapunov functionals for the QPP systems involved will be 
possible thanks to their Poisson structure, which allows the use of the 
energy-Casimir method (see [\ref{hmrw1}] and references therein) in which 
the stability analysis of a given fixed point ${\bf x}_0$ proceeds by 
defining an {\em ansatz\/} for the Lyapunov functional, which takes the form:
\begin{equation}
    \label{ecf}
    H_C({\bf x}) = H + F(C_1, \ldots ,C_k) 
\end{equation}
where $F(z_1, \ldots ,z_k)$ is a $C^2$ real function to be determined and 
$\{ C_1, \ldots , C_k \}$ is a complete set of independent Casimir 
invariants. The method amounts to the search of one suitable $F$, by imposing 
two conditions on $H_C$: (i) $H_C$ must have a critical point at ${\bf x}_0$; 
and (ii) the second derivative of $H_C$ at ${\bf x}_0$ must be either 
positive or negative definite. Once one suitable $F$ has been found, 
stability of ${\bf x}_0$ follows automatically, and the method provides us 
with a Lyapunov functional for this point. 

The structure of the article is as follows: Section 2 is devoted to the 
establishment of the main stability theorem. Different consequences of the 
result are considered in the examples of Section 3.

%\mbox{}
\pagebreak
\begin{center} 
{\sc 2. Stability of QPP Systems}
\end{center}

Let us start by recalling the following definition, valid for normed spaces 
(see [\ref{hmrw1}]):

\mbox{}

{\sc Definition 1.} [\ref{hmrw1}] \hspace{0.2cm} {\em A given steady state 
${\bf x}_0$ of a dynamical system is said to be locally stable if for 
every $\varepsilon >0$ there exists a $\delta >0$ such that: If 
$\| {\bf x}(0) - {\bf x}_0 \| < \delta$, then $\| {\bf x}(t) - {\bf x}_0 
\| < \varepsilon$ for every $t>0$.} 

\mbox{}

In what follows, stability shall denote local stability. We give now our 
main result:

\mbox{}

{\sc Theorem 2.} \hspace{0.2cm} {\em Consider a QPP system with either $m=n$ 
and $\mid B \mid >0$, or with $m>n$. If matrix D is positive or negative 
definite, then: 

\mbox{(i)} Every fixed point belonging to the interior of the 
positive orthant is stable.

\mbox{(ii)} For every fixed point belonging to the interior of the 
positive orthant there is a Lyapunov functional of the form:
\begin{equation}
    \label{fecqmp}
    H_C = H + \sum _{i=1}^n N_k \ln x_k \;\: , \: N \in \mbox{\rm ker}(K) \: ,
\end{equation}
where H is the Hamiltonian (\ref{H}).
}

\mbox{}

{\em Proof.\/} 

The proof rests strongly on the quasimonomial formalism. The unfamiliar 
reader is referred to the basic references on the subject 
[\ref{br1}-\ref{byg1},\ref{byv1}-\ref{bvb1},\ref{pym1}].

The strategy of the proof consists in reducing the problem to the 
Lotka-Volterra representative and analyzing there the stability of the fixed 
points. The resulting criteria and Lyapunov functionals are then mapped back 
into the original system. 

For the sake of clarity, we omit in what follows the proofs of the auxiliary 
lemmas, which can be found in the Appendix.

\mbox{}

{\em Proof of the first statement of Theorem 2.}

We begin by examining the behaviour of stability properties under 
embeddings. Consider an arbitrary quasipolynomial system with $m>n$:
\begin{equation}
   \dot{x}_i = x_i \left(\lambda _{i} + \sum_{j=1}^{m}A_{ij}\prod_{k=
      1}^{n}x_k^{B_{jk}} \right) , \;\:\;\: i = 1, \ldots ,n 
   \label{glv2}
\end{equation}
Let $\tilde{A}$, $\tilde{B}$ and $\tilde{\lambda}$ be the matrices of the 
expanded system which are defined in the following way: 
\begin{equation}
\label{exps}
   \tilde{A} = \left( \begin{array}{c} 
                  A_{n \times m} \\ O_{(m-n) \times m} 
               \end{array} \right)  \: , \;\;\:
   \tilde{B} = \left( \begin{array}{cc}
                  B_{m \times n} \mid B'_{m \times (m-n)}
               \end{array} \right)  \: , \;\;\:
   \tilde{\lambda} = \left( \begin{array}{c}
                  \lambda _{n \times 1} \\ O_{(m-n) \times 1}
               \end{array} \right)
\end{equation}
where we have explicitly indicated by means of indexes the sizes of the 
submatrices for the sake of clarity, $O$ denotes a null matrix, $B'$ is a 
matrix of arbitrary entries chosen in such a way that $\mid \tilde{B} \mid 
>0$, and $x_i =1$ for $i= n+1, \ldots, m$. 

\mbox{}

{\sc Lemma 1.} \hspace{0.2cm} 
{\em Let ${\bf x}_0$ $=$ $(x_{01}, \ldots ,x_{0n})$, 
with $x_{0i} >0$ for $i=1, \ldots ,n$, be a phase-space point of (\ref{glv2}), 
and let $\tilde{{\bf x}}_0$ $=$ $(x_{01}, \ldots , x_{0m})$,  with 
$x_{0i} =1$ for $i=n+1, \ldots ,m$, be the corresponding phase-space point of 
the expanded system. Then:

\mbox{(i)} ${\bf x}_0$ is a fixed point of (\ref{glv2}) if and only 
if $\tilde{{\bf x}}_0$ is a fixed point of the expanded system.

\mbox{(ii)} If ${\bf x}_0$ and $\tilde{{\bf x}}_0$ are fixed points, 
then ${\bf x}_0$ is stable if and only if $\tilde{{\bf x}}_0$ is stable.}

\mbox{}

We can now examine the effect of quasimonomial transformations (QMT's from 
now on) of the form:
\begin{equation}
\label{qmt}
   x_{i} = \prod_{k=1}^{n} y _{k}^{\Gamma_{ik}} , \;\:\: i=1, \ldots ,n \:\:
   , \;\:\: \mid \Gamma \mid > 0
\end{equation}

\mbox{}

{\sc Lemma 2.} \hspace{0.2cm} 
{\em Given a quasipolynomial system of the form 
(\ref{glv2}) with $m \geq n$, and a stable fixed point ${\bf x}_0$ belonging 
to the positive orthant, the image of ${\bf x}_0$ under an arbitrary QMT 
of the form (\ref{qmt}) is also stable.}

\mbox{}

In particular, Lemma 2 applies to the expanded QP system (\ref{exps}). Let 
us choose a QMT such that $\Gamma$ in (\ref{qmt}) is given by $\tilde{B}^{-1}$ in 
(\ref{exps}). The result is a new QP system with characteristic matrices:
\begin{equation}
\label{mlv}
    \tilde{A}' = \tilde{B} \cdot \tilde{A} \;\:\: , \;\:\:
    \tilde{B}' = I \;\:\: , \;\:\:
    \tilde{\lambda}' = \tilde{B} \cdot \tilde{\lambda} 
\end{equation}
and thus a Lotka-Volterra system ($\tilde{B}'$ is the identity 
matrix). The inverse transformation, leading from (\ref{mlv}) to (\ref{exps}) 
is also a QMT, thus validating Lemma 2 for (\ref{mlv}). 

Alternatively, in the case $m=n$ no embedding is to be performed and Lemma 2 
is applied directly to the original flow setting $\Gamma = B^{-1}$.

In either case ($m>n$ or $m=n$) we have reduced the stability problem to that 
corresponding to the Lotka-Volterra representative: If we establish stability 
for the corresponding fixed point of the Lotka-Volterra system, the steady 
state of the original flow will automatically be stable. Note that these 
considerations hold irrespectively of the fact that now the Lotka-Volterra 
representative may have an infinity of fixed points, even if this is not the 
case for the original flow. 

Let us then consider an arbitrary $m$-dimensional QPP system of 
Lotka-Volterra form. Since the tildes and primes appearing in (\ref{mlv}) 
will not be necessary in what follows, we drop them for the sake of clarity. 
We then have $A = K \cdot D$, $\lambda = K \cdot L$ and, according to 
[\ref{byv5}], rank($A$) $=$ rank($K$) $\equiv r \leq m$. Steady states are 
given in parametric form by: 
\begin{equation}
      {\bf x}_0 (N) = - D^{-1} \cdot (L-N) \;\: , \; N \in \mbox{ker}(K)
\end{equation}
We can now turn to the characterization of stability of steady-states by 
means of the energy-Casimir method. The ($m-r$) independent Casimir functions 
are of the form:
\begin{equation}
      C^{(N)} = \sum _{j=1}^m N_j \ln x_j \;\: , \; N \in \mbox{ker}(K)
\end{equation}
and we can accordingly take the following convenient form for the 
energy-Casimir functional: 
\begin{equation}
\label{fecd}
H_C \equiv \sum _{j=1}^m \left( D_{jj} x_j + (L_j + N_j) \ln x_j \right) \;\: , 
\end{equation}
where $N \in \mbox{ker}(K)$ is to be determined. Let us concentrate on a 
particular steady state ${\bf x}_0^* = -D^{-1} \cdot (L-N_0)$. We can state: 

\mbox{}

{\sc Lemma 3.} \hspace{0.2cm} 
{\em If the entries of $(L-N_0)$ are either all 
positive or all negative, then ${\bf x}_0^*$ is stable.}

\mbox{}

Now notice that:
\begin{equation}
      L-N_0 = -D \cdot {\bf x}_0^*
\end{equation}
Since we consider only steady states belonging to the positive orthant, 
Lemma 3 can be equivalently formulated in terms of positiveness or 
negativeness of matrix $D$. Since $D$ is invariant under QMTs and embeddings 
[\ref{byv5}], the same result is valid for the original QPP system and the 
first part of Theorem 2 is demonstrated.

%\mbox{}

\pagebreak
{\em Proof of the second statement of Theorem 2.}

The energy-Casimir functional (\ref{fecd}) is mapped into a functional of the 
form (\ref{fecqmp}) for the original QPP system [\ref{byv5}]. We need to 
prove, however, that (\ref{fecqmp}) is also an energy-Casimir functional. 
This is done in the following two lemmas:

\mbox{}

{\sc Lemma 4.} \hspace{0.2cm}
{\em Every QMT of the form (\ref{qmt}) maps an energy-Casimir functional into 
an energy-Casimir functional.}

\mbox{}

And finally:

\mbox{}

{\sc Lemma 5.} \hspace{0.2cm} 
{\em The property of being an energy-Casimir functional is preserved in the 
process of decoupling the $(m-n)$ variables of the embedding.}

\mbox{}

This completes the proof of Theorem 2. 

\mbox{}

{\em Remark 1.\/} The stable character of the steady state is 
independent of important features of the system, such as the degree of 
nonlinearity or the number of fixed points present in the positive orthant. 
This implies that there are certain degrees of freedom available in the 
Hamiltonian which can be varied without destroying the stability of motion. 
This has relevant consequences that we shall illustrate in the next section.
\vspace{2mm}

{\em Remark 2.\/} The criterion in Theorem 2 can be verified straightforwardly 
by simple inspection of the Hamiltonian. In particular, a precise knowledge 
of the coordinates of the fixed point(s) is not required.
\vspace{2mm}

{\em Remark 3.\/} In the specific case of conservative Lotka-Volterra 
equations, we have from (\ref{kdl}) that $B$ is the identity matrix and then 
$A = K \cdot D$. Therefore, if the hypothesis of Theorem 2 is verified then 
there exists a diagonal positive definite matrix $\bar{D}$, which is the 
absolute value of $D$, such that $\bar{D} \cdot A + A^T \cdot \bar{D} = 0$ 
due to the skew-symmetry of $K$. Accordingly, the classical stability 
criterion for conservative Lotka-Volterra systems is implied by Theorem 2 and 
now appears as a particular case.

%\mbox{}
\pagebreak
\begin{center} 
{\sc 3. Examples}
\end{center}

{\sc Example 1.} \hspace{0.2cm} 
We first consider Volterra's [\ref{v1}] predator-prey equations:
\begin{equation}
\label{lv2d}
    \begin{array}{ccl}
    \dot{x}_1 & = & x_1 ( a - bx_2) \vspace{2mm} \\
    \dot{x}_2 & = & x_2 (-d + cx_1) 
    \end{array}
\end{equation}
Here $a$, $b$, $c$ and $d$ are positive constants. This system is QPP with: 
\begin{equation}
\label{mat2d}
   K = \left( \begin{array}{cc}
         0 & 1 \\ -1 & 0 
       \end{array} \right) \; , \:\; 
   D = \left( \begin{array}{cc} 
         -c & 0 \\ 0 & -b 
       \end{array} \right) \; , \:\; 
   L = \left( \begin{array}{c} 
         d \\ a 
       \end{array} \right) 
\end{equation}
The Hamiltonian is:
\begin{equation}
   \label{h1e1}
   H(x_1,x_2) = -cx_1 -bx_2 +d \ln x_1 +a \ln x_2
\end{equation}
It is well known that there is a unique fixed point in the positive orthant, 
which is stable. We can immediately verify this from the point of view of 
Theorem 2, since $D$ in (\ref{mat2d}) is negative definite. Therefore the 
steady state is stable. Moreover, (\ref{h1e1}) is a Lyapunov functional for 
it, since flow (\ref{lv2d}) is symplectic.

\mbox{}

{\sc Example 2.} \hspace{0.2cm} Taking the system of Example 1 as starting 
point, let us now consider the following generalization of the Hamiltonian: 
\begin{equation}
   \label{h2e1}
   H(x_1,x_2) = -cx_1^{\alpha}x_2^{\beta} -bx_1^{\gamma}x_2^{\delta} + 
                 d \ln x_1 + a \ln x_2
\end{equation}
Now the equations become 
\begin{equation}
\label{nlg1}
          \begin{array}{ccl} 
      \dot{x}_1 & = & x_1 ( a - \beta c x_1^{\alpha}x_2^{\beta} - \delta b
                           x_1^{\gamma}x_2^{\delta}) \vspace{2mm} \\
      \dot{x}_2 & = & x_2 (-d + \alpha c x_1^{\alpha}x_2^{\beta} + \gamma b
                           x_1^{\gamma}x_2^{\delta}) \vspace{2mm}
          \end{array}
\end{equation}

Let us assume that $\alpha$, $\beta$, $\gamma$ and $\delta$ are all positive. 
Since in Volterra's model $\alpha$ and $\delta$ are greater than $\beta$ and  
$\gamma$, we shall also extend this requirement here and consider:
\begin{equation}
   \mid B \mid = \alpha \delta - \beta \gamma > 0
\end{equation}
Within these assumptions, which are not very restrictive, it is not difficult 
to prove that there exists a unique fixed point inside the positive orthant 
if and only if:
\begin{equation}
   \label{1.i}
     \frac{\delta}{\gamma} > \frac{a}{d} > \frac{\beta}{\alpha}
\end{equation}
We have that matrix $D$ retains the same form than in (\ref{mat2d}). 
Therefore, according to Theorem 2 the point is stable, Hamiltonian 
(\ref{h2e1}) is also a Lyapunov functional of the generalized system (given 
that (\ref{nlg1}) is a symplectic flow) and (\ref{1.i}) remains as the only 
condition both for the positiveness of the fixed point and for its stability. 

It is clear that the generalized Hamiltonian (\ref{h2e1}) must 
incorporate dynamical features not present in Volterra's model. To see this, 
we first put (\ref{nlg1}) into classical Hamiltonian form by means of 
transformation $y_i = \ln x_i$, for $i=1,2$ (see [\ref{byv5}] for the general 
reduction algorithm of QPP systems into the Darboux canonical form). 
After that, we perform a phase-space translation with the new axes centered 
in the steady state: $y_i = y_i^0 + \varepsilon _i$, $i=1,2$. Finally, we 
consider the case of small oscillations around the steady state and neglect 
terms of order $\varepsilon ^3$. The resulting system has the following 
Hamiltonian: 
\begin{equation}
\label{h3e1}
    H(\varepsilon _1,\varepsilon _2) = \mu _1 \varepsilon _1^2 + \mu _2 \varepsilon _2^2 
                                 + 2 \mu \varepsilon _1 \varepsilon _2 \;\: ,
\end{equation}
where $\mu_1$, $\mu_2$ and $\mu$ are negative constants. 

We shall first consider a particular case of (\ref{nlg1}) 
\begin{equation}
\label{pc1}
   \begin{array}{ccl}
    \dot{x}_1 & = & x_1 ( a - (1+ \delta ^*) b x_2^{1+\delta ^*}) \vspace{2mm} \\
    \dot{x}_2 & = & x_2 (-d + (1+ \alpha ^*) c x_1^{1+\alpha ^*})
   \end{array}
\end{equation}
where $\alpha ^*$ and $\delta ^*$ are greater than $-1$. It is a simple task to 
demonstrate that for (\ref{pc1}) $\mu =0$ in 
(\ref{h3e1}), and then the trajectories are ellipses aligned with the 
coordinate axes, similarly to what occurs in Volterra's case. However, the 
frecuency of the oscillations is now generalized to:
\begin{equation}
      \omega = \sqrt{(1+ \alpha ^*)(1+ \delta ^*)ad}
\end{equation}
If $\alpha ^*$ and $\delta ^*$ remain small, $\omega$ is of the order of 
Volterra's frecuency $\omega _0 = \sqrt{ad}$. In the most general case, 
$\omega$ can take any positive value, and is not restricted to any particular 
range. 

There are some additional features not present in Volterra's model which are 
due to the off-diagonal terms in matrix $B$. These are related to the phase 
shift between the oscillations of the predator and the prey. To see this, let 
us turn back to the general Hamiltonian (\ref{h3e1}) for the case of small 
oscillations. It is well known that there exists a canonical transformation, 
which is a rotation of angle $\phi$ of the axes, such that in the new 
variables the Hamiltonian is 
\begin{equation}
   H(\xi _1, \xi _2) = \lambda _1 \xi _1^2 + \lambda _2 
   \xi _2^2 \;\: ,
\end{equation}
where $\lambda _1 $ and $\lambda _2$ are the eigenvalues of the ellipse. The 
solution for $(\xi _1, \xi _2)$ is straightforward. Then, if we 
transform back into the variables $(\varepsilon _1 , \varepsilon _2)$, a 
simple calculation shows that the phase shift between the predator and the 
prey is just:
\begin{equation}
   \Phi (\rho,\phi) = \frac{\pi}{2} + \arctan ( \rho \tan \phi ) -
    \arctan ( \rho ^{-1} \tan \phi )
\end{equation}
where $\rho = \sqrt{\lambda _1 / \lambda _2}$. Thus, we now have phase shifts 
which may be different to $\pi /2$, which is the classical Volterra value 
($\phi = 0$). Notice that, in the neighbourhood of $\phi = 0$ we have:
\begin{equation}
   \Phi (\rho,\phi) = \frac{\pi}{2} + \left( \rho - \frac{1}{ \rho } \right) 
                      \phi + o( \phi ^3)
\end{equation}
Therefore, if the eigenvalues do not have exactly the same magnitude (which 
is a reasonable assumption) these models can reproduce, in particular, 
a whole range of phase shifts centered around $\pi /2$. This is consistent 
with observed time series in predator-prey systems (see, for example, 
[\ref{hass1}, pp. 60, 92] and [\ref{mur1}, p. 67]) in which the average 
phase shifts may differ from $\pi /2$. 

We can then conclude that generalization (\ref{h2e1}) accounts for additional 
features observed in real systems, while retaining the advantages and the 
basic framework provided by a Hamiltonian formulation. 

\mbox{}

{\sc Example 3.} \hspace{0.2cm} 
We shall start again with the Lotka-Volterra equations (\ref{lv2d}). Let us 
now consider the addition to the Hamiltonian (\ref{h1e1}) of two extra 
nonlinear terms: 
\begin{equation}
      H(x_1,x_2) = -cx_1 -bx_2 + \sigma _1 x_1^{\alpha} + 
      \sigma _2 x_2^{\beta} + d \ln x_1 +a \ln x_2 
\end{equation}
with both $\alpha$ and $\beta$ positive and different from 1. Notice that 
matrix $D$ is:
\begin{equation}
   D = \left( \begin{array}{cccc} 
         -c & 0 & 0 & 0 \\ 0 & -b & 0 & 0 \\ 
         0 & 0 & \sigma _1 & 0 \\ 0 & 0 & 0 & \sigma _2
       \end{array} \right)
\end{equation}
The resulting generalized equations are: 
\begin{equation}
\label{e2eq1}   
   \begin{array}{ccl}
   \dot{x}_1 & = & x_1 ( a - bx_2 + \beta \sigma_2 x_2^{\beta}) \vspace{2mm} \\
   \dot{x}_2 & = & x_2 (-d + cx_1 - \alpha \sigma_1 x_1^{\alpha}) 
   \end{array}
\end{equation}
Before considering the existence of steady states, note from the form of 
$H$ and $D$ and from Theorem 2 that every fixed point of the positive 
orthant is stable if $\sigma _1 <0$ and $\sigma _2 <0$, independently of the 
values of $\alpha$ and $\beta$. Let us assume that this is the case. It is 
then simple to prove that there exists a unique point in the interior of the 
positive orthant which verifies the fixed point conditions:
\begin{equation}
   \begin{array}{ccl}
            cx_1 - \alpha \sigma_1 x_1^{\alpha} & = &  d    \\
            bx_2 -  \beta \sigma_2 x_2^{\beta}  & = &  a 
   \end{array}
\end{equation}
Therefore there is a unique steady state inside the positive orthant, it 
is stable and $H$ is a Lyapunov functional for it. The analytic 
determination of the coordinates of the point may be a nontrivial problem, 
since $\alpha$ and $\beta$ are real constants in general. However, it is now 
possible to establish stability even without knowing the exact position of 
the point, but only by demonstrating its existence. 
         
\mbox{}

{\sc Example 4.} \hspace{0.2cm} 
We shall finally look upon the following system, characterized by Nutku 
[\ref{nut1}]:
\begin{equation}
\label{lv3d}
   \begin{array}{ccl}
   \dot{x}_1 &=& x_1( \rho + cx_2 +  x_3) \vspace{2mm}  \\
   \dot{x}_2 &=& x_2( \mu  +  x_1 + ax_3) \vspace{2mm}  \\
   \dot{x}_3 &=& x_3( \nu  + bx_1 +  x_2) \vspace{2mm}  
   \end{array}
\end{equation}
As Nutku has pointed out, this is a Poisson system if 
\begin{equation}
\label{clv3d}
   abc=-1 \; , \;\;\; \nu = \mu b - \rho ab
\end{equation}
In fact, if conditions (\ref{clv3d}) hold the system is QPP with Hamiltonian:
\begin{equation}
\label{hlv3d}
   H = ab x_1 + x_2 - a x_3 + \nu \ln x_2 - \mu \ln x_3
\end{equation}
The associated QPP matrices are:
\begin{equation}
\label{mlv3d}
   K = \left( \begin{array}{ccc}
          0 &  c & bc \\
         -c &  0 & -1 \\
        -bc &  1 & 0 
       \end{array} \right)
\; , \;\;\; 
   D = \left( \begin{array}{ccc}
         ab &  0 &  0 \\
         0  &  1 &  0 \\
         0  &  0 & -a
       \end{array} \right)
\; , \;\;\; 
   L = \left( \begin{array}{c}
         0 \\ \nu \\ - \mu 
       \end{array} \right)
\end{equation}
System (\ref{lv3d}), being odd-dimensional, falls out of the scope of the 
traditional Hamiltonian domain. However, the previous results hold in this 
context as well. If we apply Theorem 2 to this case, from $D$ in (\ref{mlv3d}) 
together with (\ref{clv3d}) we immediately obtain that the fixed points of 
the positive orthant are stable if $a<0$, $b<0$ and $c<0$. Notice that system 
(\ref{lv3d}) has an infinite number of fixed points, so stability is 
simultaneously demonstrated for all those belonging to int$\{I \!\! R^3_+\}$. 

Notice also that the flow is not symplectic, and we have one independent 
Casimir invariant:
\begin{equation}
     C = ab \ln x_1 -b \ln x_2 + \ln x_3 = \mbox{constant} 
\end{equation}
Thus, according to (\ref{fecqmp}) the Lyapunov functional of every positive 
steady state will be of the form:
\begin{equation}
   H_C = H + \kappa C = ab x_1 + x_2 - a x_3 + 
         \kappa ab \ln x_1 + ( \nu - \kappa b) \ln x_2 + 
         ( \kappa - \mu ) \ln x_3 \; , \:\; \kappa \in I \!\! R
\end{equation}
Obviously, the Lyapunov functional (i.e., the appropriate value of the 
parameter $\kappa$) will be different for every fixed point and can be 
determined without difficulty by following the constructive procedure given 
in the proof of Theorem 2. We do not elaborate further on this issue for the 
sake of conciseness. 

Finally, notice that the flow can be easily generalized to account for higher 
order nonlinearities while preserving stability, by means of the same 
techniques employed in Examples 2 and 3. Such techniques are completely 
general.

\pagebreak
\begin{center}
{\sc Appendix}
\end{center}

{\em Proof of Lemma 1.} \hspace{0.2cm}
For (i), we have:
\[
   \sum_{j=1}^m \tilde{A}_{ij} \prod _{k=1}^m (x_{0k})^{\tilde{B}_{jk}} 
   + \tilde{\lambda}_i = 0 \;\: , \: \forall i = 1, \ldots ,m \Longrightarrow 
\]
\begin{equation}
   \sum_{j=1}^m A_{ij} \prod _{k=1}^n (x_{0k})^{B_{jk}} + \lambda_i = 0  
   \;\: , \: \forall i 
\end{equation}
The converse follows after noting that the sense of these implications can 
be reversed.

The proof of (ii) is a consequence of the fact that the removal or addition 
of variables of constant value 1 does not affect the stable character of the 
point. {\bf Q.E.D.}

\mbox{}

{\em Proof of Lemma 2.} \hspace{0.2cm}
It is a consequence of the fact that QMTs (\ref{qmt}) are 
orientation-preserving diffeomorphisms and therefore relate topologically 
orbital equivalent systems. {\bf Q.E.D.}

\mbox{}

{\em Proof of Lemma 3.} \hspace{0.2cm}
The gradient of the energy-Casimir functional vanishes identically at 
${\bf x}_0^*$ if we set $N = -N_0$ in $H_C$. For the second part of the 
criterion, we note that the Hessian of $H_C$ at ${\bf x}_0^*$ is diagonal 
due to the simple form of $H$ in the case of Lotka-Volterra equations, and 
takes the value 
\begin{equation}
   \mbox{Hess}(H_C \mid _{{\bf x}_0^*}) = \mbox{diag} \left(
         \frac{(N_0-L)_1}{(x_{01}^*)^2}, \ldots , 
         \frac{(N_0-L)_m}{(x_{0m}^*)^2} \right)
\end{equation}
{\bf Q.E.D.}

\mbox{}

{\em Proof of Lemma 4.} \hspace{0.2cm}
It is a simple consequence of the chain rule for $C^2$ functions of $m$ real 
arguments. {\bf Q.E.D.}

\mbox{}

{\em Proof of Lemma 5.} \hspace{0.2cm}
Clearly, if the gradient of $H_C$ vanishes at $\tilde{{\bf x}}_0$, the 
gradient of the $n$-dimensio\-nal restriction of $H_C$ will also vanish at 
${\bf x}_0$. Similarly, the Hessian of the restriction of $H_C$ will be 
a $n \times n$ minor of the Hessian of $H_C$, corresponding to the first 
$n$ rows and columns. Consequently, the Hessian of the restriction will also 
be definite. {\bf Q.E.D.}

\pagebreak
\begin{center}
{\sc References}
\end{center}
\begin{enumerate}
\item L. Brenig, Complete factorisation and analytic solutions of generalized 
      Lotka-Volterra equations, {\em Phys. Lett. A\/} 133 (1988), 378-382.
      \label{br1}
\item L. Brenig, and A. Goriely, Universal canonical forms for time-continuous 
      dynamical systems, {\em Phys. Rev. A\/} 40 (1989), 4119-4122.
      \label{byg1}
\item A. Gasull, A. Guillamon, C. Li and Z. Zhang, Study of Perturbed 
      Lotka-Volterra Systems Via Abelian Integrals, {\em J. Math. Anal. 
      Appl.\/} 198 (1996), 703-728.
      \label{gas1}
\item N. S. Goel, S. C. Maitra and E. W. Montroll, On the Volterra and Other 
      Nonlinear Models of Interacting Populations, {\em Rev. Mod. Phys.\/}
      43 (1971), 231-276.
      \label{gmm1}
\item B. S. Goh, Global stability in many-species systems, {\em Amer. 
      Naturalist\/} 111 (1977), 135-143.
      \label{goh1}
\item B. S. Goh, Sector Stability of a Complex Ecosystem Model, {\em Math. 
      Biosci.\/} 40 (1978), 157-166.
      \label{goh2}
\item H. Goldstein, ``Classical Mechanics,'' Second Edition, Addison-Wesley,
      Reading (Massachusetts), 1962.
      \label{gold}
\item G. W. Harrison, Global stability of food chains, {\em Amer. 
      Naturalist\/} 114 (1979), 455-457.\label{har1}
\item M. P. Hassell, ``The Dynamics of Arthropod Predator-Prey Systems,'' 
      Princeton University Press, Princeton (New Jersey), 1978.
      \label{hass1}
\item B. Hern\'{a}ndez--Bermejo and V. Fair\'{e}n, Nonpolynomial vector 
      fields under the Lotka-Volterra normal form, {\em Phys. Lett. A\/} 
      206 (1995), 31-37.
      \label{byv1}
\item B. Hern\'{a}ndez--Bermejo and V. Fair\'{e}n, Lotka-Volterra 
      Representation of General Nonlinear Systems, {\em Math. Biosci.\/} 140
      (1997), 1-32.
      \label{byv2}
\item B. Hern\'{a}ndez--Bermejo and V. Fair\'{e}n, Hamiltonian structure and 
      Darboux theorem for families of generalized Lotka-Volterra systems,
      {\em J. Math. Phys.\/} 39 (1998), 6162-6174.
      \label{byv5}
\item B. Hern\'{a}ndez--Bermejo, V. Fair\'{e}n and L. Brenig, Algebraic 
      recasting of nonlinear systems of ODEs into universal formats, {\em J. 
      Phys. A, Math. Gen.\/} 31 (1998), 2415-2430.\label{bvb1}
\item J. Hofbauer and K. Sigmund, ``The Theory of Evolution and Dynamical 
      Systems,'' Cambridge University Press, Cambridge, 1988.
      \label{hs1}
\item D. D. Holm, J. E. Marsden, T. Ratiu and A. Weinstein, Nonlinear 
      Stability of Fluid and Plasma Equilibria, {\em Phys. Rep.\/} 123 
      (1985), 1-116.
      \label{hmrw1}
\item S. B. Hsu, The Application of the Poicar\'{e}-Transform to the 
      Lotka-Volterra Model, {\em J. Math. Biol.\/} 6 (1978), 67-73.
      \label{hsu1}
\item E. H. Kerner, A statistical mechanics of interacting species, {\em Bull. 
      Math. Biophys.\/} 19 (1957), 121-146.
      \label{ker1}
\item E. H. Kerner, ``Gibbs Ensemble: Biological Ensemble,'' Gordon and 
      Breach, New York, 1972.
      \label{ker2}
\item E. H. Kerner, Note on Hamiltonian format of Lotka-Volterra dynamics, 
      {\em Phys. Lett. A\/} 151 (1990), 401-402.
      \label{ker3}
\item E. H. Kerner, Comment on Hamiltonian structures for the $n$-dimensional 
      Lotka-Volterra equations, {\em J. Math. Phys.\/} 38 (1997), 1218-1223.
      \label{ker4}
\item N. Krikorian, The Volterra model for three-species predator-prey 
      systems: Boundedness and stability, {\em J. Math. Biol.\/} 7 (1979), 
      117-132.
      \label{kri1}
\item J. P. LaSalle, Some Extensions of Liapunov's Second Method, {\em IRE 
      Trans. Circuit Theory\/} 7 (1960), 520-527.
      \label{las0}
\item J. P. LaSalle, Stability Theory for Ordinary Differential Equations, 
      {\em J. Differ. Equations\/} 4 (1968), 57-65.
      \label{las1}
\item X.-Z. Li, C.-L. Tang and X.-H. Ji, The Criteria for Globally Stable 
      Equilibrium in $n$-Dimensional Lotka-Volterra Systems, {\em J. Math. 
      Anal. Appl.\/} 240 (1999), 600-606.\label{ltj1}
\item D. O. Logofet, ``Matrices and Graphs. Stability Problems in Mathematical
      Ecology,'' CRC Press, Boca Raton (Florida), 1993.
      \label{log1}
\item A. J. Lotka, ``Elements of Mathematical Biology,'' Dover, New York, 
      1956.
      \label{lot1}
\item Z.-Y. Lu and Y. Takeuchi, Global Stability Conditions for 
      Three-Dimensional Lotka-Volterra Systems, {\em Appl. Math. Lett.\/} 7
      (1994), 67-71.
      \label{lyt1}
\item Z. Lu and Y. Takeuchi, Global Dynamic Behavior for Lotka-Volterra 
      Systems with a Reducible Interaction Matrix, {\em J. Math. Anal. 
      Appl.\/} 193 (1995), 559-572.
      \label{lyt2}
\item J. D. Murray, ``Mathematical Biology,'' Second Edition, Springer-Verlag, 
      Berlin, 1993.\label{mur1}
\item Y. Nutku, Hamiltonian structure of the Lotka-Volterra equations, 
      {\em Phys. Lett. A\/} 145 (1990), 27-28.
      \label{nut1}
\item P. J. Olver, ``Applications of Lie Groups to Differential Equations,'' 
      Second Edition, Springer-Verlag, New York, 1993.
      \label{olv1}
\item M. Peschel and W. Mende, ``The Predator--Prey Model. Do we live in a 
      Volterra World?'' Springer-Verlag, Vienna, 1986.
      \label{pym1}
\item M. Plank, Hamiltonian structures for the $n$-dimensional 
      Lotka-Volterra equations, {\em J. Math. Phys.\/} 36 (1995), 3520-3534.
      \label{pla1}
\item M. Plank, On the Dynamics of Lotka-Volterra Equations Having an 
      Invariant Hyperplane, {\em SIAM J. Appl. Math.\/} 59 (1999), 1540-1551.
      \label{pla2}
\item R. Redheffer, A new class of Volterra differential equations for which 
      the solutions are globally asymptotically stable, {\em J. Differ. 
      Equations\/} 82 (1989), 251-268.
      \label{red1}
\item Y. Takeuchi, ``Global Dynamical Properties of Lotka-Volterra Systems,'' 
      World Scientific, Singapore, 1996.
      \label{tak1}
\item Y. Takeuchi and N. Adachi, The Existence of Globally Stable Equilibria 
      of Ecosystems of the Generalized Volterra Type, {\em J. Math. Biol.\/}
      10 (1980), 401-415.
      \label{ta1}
\item Y. Takeuchi and N. Adachi, Stable Equilibrium of Systems of Generalized 
      Volterra Type, {\em J. Math. Anal. Appl.\/} 88 (1982), 157-169.
      \label{ta2}
\item Y. Takeuchi, N. Adachi and H. Tokumaru, The Stability of Generalized 
      Volterra Equations, {\em J. Math. Anal. Appl.\/} 62 (1978), 453-473.
      \label{tat1}
\item Y. Takeuchi, N. Adachi and H. Tokumaru, Global Stability of Ecosystems 
      of the Generalized Volterra Type, {\em Math. Biosci.\/} 42 (1978), 
      119-136.
      \label{tat2}
\item V. Volterra, ``Le\c{c}ons sur la Th\'{e}orie Math\'{e}matique de la 
      Lutte pour la Vie,'' Gauthier Villars, Paris, 1931.
      \label{v1}
\end{enumerate}

\end{document}